\documentclass[10pt,a4paper]{article}

\usepackage{graphicx}
\usepackage{amsmath,latexsym}
\usepackage{amssymb}
\usepackage{amsfonts}
\usepackage{amsthm}

\newtheorem{thm}[subsection]{Theorem}
\newtheorem{lem}[subsection]{Lemma}
\newtheorem{prop}[subsection]{Proposition}
\newtheorem{cor}[subsection]{Corollary}

\newtheorem{defn}[subsection]{Definition}

\newtheorem{examples}[subsection]{Examples}

\newcommand{\Cros}{\begin{picture}(30,15) \put(0,15){\line(2,-1){30}}
       \multiput(0,0)(20,10){2}{\line(2,1){10}} \end{picture}}
\newcommand{\cros}{\begin{picture}(30,15) \put(0,15){\line(2,-1){30}}
       \multiput(0,0)(20,10){2}{\line(2,1){10}} \end{picture}}

\title{Centralizer of braids and Fibonacci numbers\thanks{This research is
supported by the grant of Higher Education Commission, Pakistan, Ref
no: PM-IPFP/HRD/HEC/2010/1508.\newline \textbf{Keywords:}  simple
braids, Fibonacci numbers, commuting graphs.\newline 2000
\textbf{Mathematics Subject Classification:} 11B39, 05A15,
05A05.\newline \textbf{Corresponding author:} uali@bzu.edu.pk}}
\author{Usman Ali, Azeem Haider}
\date{}
\begin{document}
\maketitle
\vspace{-7mm}
\begin{center}
{\it\small Center for Advanced Studies in Pure and Applied
Mathematics, Bahauddin Zakaria University, Multan, Pakistan. {\it
E-mail:} \{uali, azeemhaider\}@bzu.edu.pk}
\end{center}
%%%%%%%%%%%%%%%%%%%%%%%%%%%%%%%%
%%%%%%%%%%%%%%%%%%%%%%%%%%%%%%%%
%%%%%%%%%%%%%%%%%%%%%%%%%%%%%%%%
\begin{abstract}
The paper encloses computation of simple centralizer of simple
braids and their connection with Fibonacci numbers. Planarity of
some commuting graphs is also discussed in the last section.
\end{abstract}
\section{\textbf{Introduction}} The set of
\textit{positive n-braids} $\mathcal{MB}_{n}$ is defined by the
classical presentation\cite{art}:
\begin{equation}\label{classical}
\mathcal{MB}_{n}=\left\langle x_{1},x_{2},\cdots,x_{n-1}:
   \begin{array}{l}
\,\,x_{i+1}\,x_{i}\,x_{i+1}\,\,=\,x_{i}\,x_{i+1}\,x_{i}\\
\,\,x_{i}\,x_{j}\,=\,x_{j}\,x_{i}\,\,$for$\,\,|i-j|\geq 2
\end{array}\right\rangle.
 \end{equation}
 In other words, a positive braid $\alpha$ is a word in the set of
 generators $\{x_1,x_2,\ldots, x_{n-1}\}:$
\begin{center}\begin{picture}(365,30)                           %% fig 5a  %%
\put(95,30){$ 1 $} \put(80,17){$ x_i $} \put(122,30){$2$}
\put(157,30){$ i $} \put(180,30){$ i+1 $} \put(212,30){$ n-1 $}
\put(250,30){$ n $} \put(98,07){\line (0,1){20}}
\put(253,07){\line(0,1){20}} \put(132,15){$\ldots$}
\put(200,15){$\ldots$} \put(160,10){\Cros }
\put(125,07){\line(0,1){20}} \put(225,07){\line(0,1){20}}
\end{picture}\end{center}

\noindent such that $\alpha\in\{x^{r_1}_{i_1}\cdots
x^{r_k}_{i_k}:r_1,\ldots,r_k\geq 0\}$ and two words are equivalent
if they can be transferred into each other by finitely many
applications of relations given in~(\ref{classical}), for example
$x_4x_3^2x_2x_3$ and $x_2x_4x_3x_2^2$ are equivalent in
$\mathcal{MB}_{5}.$ The braid monoids $\mathcal{MB}_{n}$ are
embedded in their corresponding braid groups $\mathcal{B}_{n}$,
defined by the same presentation~(\ref{classical}) (see~\cite{gar}).
A \textit{simple braid} contains a letter $x_i$ at most once
(see~\cite{rb}) and
 was shown in~\cite{rbr} that the number of simple braids in $\mathcal{SB}_{n}$
 is Fibonacci number $F_{2n-1},$ where
$(F_0,F_1,F_2,F_3,F_4,F_5\ldots)=(0,1,1,2,3,5,\ldots).$

 Three kinds of divisors of $\beta\in\mathcal{MB}_{n}$ are defined: divisors of
 $\beta$ $(\gamma|\beta),$ $\textrm{Div}(\beta)=\{\gamma \in \mathcal{MB}_{n}:\, \textrm{there
 exist}\,
 \delta,\varepsilon \in\mathcal{MB}_{n}, \beta=\delta\gamma\varepsilon
 \}$; left divisors of $\beta$ $(\gamma|_L\beta),$ $\textrm{Div}_L(\beta)=\{\gamma \in \mathcal{MB}_{n}:\,\textrm{there
 exists }\,
 \varepsilon\in\mathcal{MB}_{n}, \beta=\gamma\varepsilon
 \}$; and right divisors of $\beta$ $(\gamma|_R\beta),$ $\textrm{Div}_R(\beta)= \{\gamma \in \mathcal{MB}_{n}:\, \textrm{there
 exists }\,
 \delta \in\mathcal{MB}_{n},$ $\beta=\delta\gamma\}.$ Clearly
$\textrm{Div}_L(\beta)\cup\textrm{Div}_R(\beta)\subseteq\textrm{Div}(\beta).$

The set $\mathcal{SB}_{n}$ is a proper subset of
$\textrm{Div}(\Delta_n)$, where
 $$\Delta_n=x_1(x_2x_1)\cdots(x_{n-1}x_{n-2}\cdots x_2x_1)$$ is the \textit{Garside braid}
 (see~\cite{rb} for more details). The braid $x_1x_3x_2x_4$ is simple
while the braid $x_1x_3x_2x_3$ is a non-simple divisor of
$\Delta_n$.
\begin{center}
\protect  \begin{picture}(700,80)                %%    fig 5b  %%
\put(82,62){\line(0,1){16}} \put(54,48){\line(0,1){13.3}}
\put(24,46){\line(0,1){15}}\put(24,31){\line(0,1){15}}
\put(54,18){\line(0,1){15}}
 \put(24,16){\line(0,1){15}}\put(113,33){\line(0,1){14}}
 \put(85,18){\line(0,1){15}}
 \put(143.5,32.5){\line(0,1){45}}
 \put(24,61){\Cros }
\put(83,47){\Cros} \put(54.5,33){\cros }\put(113.5,17.5){\cros }
 \put(113,62){\line(0,1){16}} \put(17,80){ $1$}
 \put(82,80){$3$}\put(82,9){$3$}\put(17,8){ $1$}
 \put(107,80){ $ 4 $ }
 \put(107,9){ $ 4 $ }
 \put(137,80){ $ 5 $ }
 \put(137,8){ $ 5 $ }
 \put(45,-20){ $ \beta=x_1x_3x_2x_4$ }
\put(48,80){ $ 2 $ } \put(48,9){ $ 2 $ }\protect
\end{picture}
\end{center}
\begin{center}
\protect  \begin{picture}(130,-80)                %%    fig 5b  %%
\put(144,79){\line(0,1){16}}
\put(174,79){\line(0,1){16}}\put(114,64){\line(0,1){14}}
\put(83,63){\line(0,1){15}}\put(83,48){\line(0,1){15}}
\put(113,35){\line(0,1){14}}
 \put(83,33){\line(0,1){15}}
 \put(174,49.5){\line(0,1){15}}
 \put(200,33){\line(0,1){63}}
 \put(83.6,78){\Cros }
\put(144,64){\Cros} \put(114,49){\cros }\put(143.8,34.5){\cros }
  \put(76,97){ $1$}
 \put(142,97){$3$}\put(142,26){$3$}\put(76,25){ $1$}
 \put(167,97){ $ 4 $ }
 \put(167,26){ $ 4 $ }
  \put(194,98){ $ 5 $ }
 \put(194,23){ $ 5 $ }
 \put(105,1){ $ \alpha=x_1x_3x_2x_3 $ }

\put(108,97){ $ 2 $ } \put(108,26){ $ 2 $ }\protect
\end{picture}
\end{center}
 \vspace{1cm}
 There is a canonical group homomorphism onto the symmetric group\linebreak $\pi:\mathcal{B}_{n}\rightarrow \Sigma_n$; for example
 $$\pi(\alpha)=\left(
     \begin{array}{ccccc}
       1 & 2 & 3 & 4& 5\\
       4& 1 & 3& 2&5\\
     \end{array}
   \right)
 .$$
 The restriction of $\pi$ to $\textrm{Div}(\Delta_n)$ is a bijection (see \cite{wpg}, chapter 9).

 The symmetric group $\Sigma_n$ also admits the following presentation in the Coxeter
 generators:
 \begin{equation}\label{symmetric}
\Sigma_{n}=\left\langle s_{1},s_{2},\cdots,s_{n-1}:
   \begin{array}{l}
\,\,s^2_i=1\\
\,\,s_{i+1}\,s_{i}\,s_{i+1}\,\,=\,s_{i}\,s_{i+1}\,s_{i}\\
\,\,s_{i}\,s_{j}\,=\,s_{j}\,s_{i}\,\,$for$\,\,|i-j|\geq 2
\end{array}\right\rangle.
 \end{equation}
The homomorphism $\pi$ can be defined by $\pi(x_i)=s_i$. The
image\linebreak $\pi(\mathcal{SB}_{n})=S\Sigma_n$ is called the set
of simple permutations in $\Sigma_n$ (see \cite{rab}).

 The braid monoids $\mathcal{MB}_{n}$ satisfy left and right cancellation laws, i.e,
$\beta\gamma=\beta\delta$ $\Rightarrow$ $\gamma=\delta$ and
$\gamma\beta=\delta\beta$ $\Rightarrow$ $\gamma=\delta$ (see
\cite{gar}). The monoid $\mathcal{MB}_{n}$ is embedded in
$\mathcal{MB}_{n+1}$ and consequently, $\alpha, \beta\in
\mathcal{MB}_{n}$ commutes in $\mathcal{MB}_{n+1}$ if and only if
$\alpha$ and $\beta$ commutes in $\mathcal{MB}_{n}.$

\begin{defn}
The \textit{simple centralizer} of $\beta\in \mathcal{SB}_{n}$ is
the set\linebreak
$C_n(\beta)=\{\gamma\in\mathcal{SB}_{n}:\beta\gamma=\gamma\beta\}$,
i.e, the intersection of centralizer of $\beta$ in $\mathcal{B}_{n}$
with $\mathcal{SB}_{n}.$ The cardinality of $C_n(\beta)$ is denoted
by $c_n(\beta).$
\end{defn}
Our first result completely describes the structure of  $C_n(x_i)$
for\linebreak $1\leq i\leq n-1$.
\begin{thm}\label{gen} For any $x_i\in \mathcal{SB}_{n},$ we have
$$C_n(x_i)=\{\beta\in\mathcal{SB}_{n}: x_j\nmid\beta\, \emph{for}\, \mid j-i\mid=1\}.$$
\end{thm}
For $x_{n-1}|\beta$ (where $n-1=$max$\{i:x_i\mid\beta\}$), our
second result describes the structure of $C_{n+m}(\beta)$ in terms
of $C_n(\beta)$.
\begin{thm}\label{struc}
If $\beta\in\mathcal{SB}_{n}$ and $x_{n-1}|\beta,$ then $\gamma\in
C_{n+m}(\beta)$ if and only if $\gamma=\gamma_1\gamma_2$
 where $m\geq 1$, $\gamma_1\in C_{n}(\beta)$, and $\gamma_2\in
E_n=\{\alpha\in \mathcal{SB}_{n+m}:x_j\nmid\alpha \emph{ for all }
j\leq n\}.$
\end{thm}
It is clear that the sets $E_n$ and $\mathcal{SB}_{m}$ have the same
number of elements and consequently the following holds.

\begin{cor}\label{cent}
If $\beta\in\mathcal{SB}_{n}$ and $x_{n-1}|\beta,$ then
$c_{n+m}(\beta)=c_n(\beta)F_{2m-1}$\linebreak for all $m \geq 1$.
\end{cor}
%In order to find the set $C_n(\beta)$ for an arbitrary $n$, we have
%to know only $C_n(\beta)$ such that $x_{n-1}|\beta.$

The value of $c_{n+m}(x_{n-1})$ is given by the following
proposition.
\begin{prop}\label{c(x)}
a) $c_{2+m}(x_1)=2F_{2m-1}$ and\\
 b) $c_{n+m}(x_{n-1})=2F_{2n-5}.F_{2m-1}$ $\forall\,\, n\geq 3 $.
\end{prop}

%Using the equality $C_n(e)=\mathcal{SB}_{n}$, we give a new proof of
%the following theorem
%
%\begin{thm}$($\cite{rbr}, Theorem~1.4.$)$ \label{nsim}
% The number of simple braids in $\mathcal{SB}_{n}$ is $F_{2n-1}.$
%\end{thm}

In the next section proofs of Theorem~\ref{gen},
Theorem~\ref{struc}, and Proposition~\ref{c(x)} are given.

An algorithm is given in~\cite{mens} to compute the centralizer of
an arbitrary braid $\beta\in\mathcal{B}_{n}.$ In the same paper,
simple elements are defined as divisors of the Garside braid
$\Delta_n$~(see definition 4 in \cite{mens}).
\begin{defn}(\cite{doug}, page 248)
 A graph is said to be \textit{planar} if it can be embedded in a plane. Otherwise it is \textit{non-planar}.
 \end{defn}
 \begin{defn}(\cite{bates1,bates2,com})
  A \textit{commuting graph} $\Gamma(H)$ associated to a group $G$ and a finite subset $H$ of
$G$ is a graph whose vertices are the elements of $H\setminus\{e\}$
and there is an edge between $g$ and $h$ if and only if $g\neq h$
and $gh=hg$.
 \end{defn}

In the last section, we discussed the planarity of three graphs:
$\Gamma(\mathcal{SB}_{n})$, $\Gamma(S\Sigma_n)$, and
$\Gamma(\Sigma_n).$ The graph $\Gamma(\Sigma_n)$ was analyzed
in~\cite{com}.

\begin{prop}\label{planner}
 a) $\Gamma(\mathcal{SB}_{n})$ is planar if and only
if $n\leq5$ and\\
b) $\Gamma(S\Sigma_n)$ and $\Gamma(\Sigma_n)$ are planar if and only
if $n\leq4$.
\end{prop}
The graph $\Gamma(\mathcal{SB}_{n})$ is a proper subgraph of
$\Gamma(S\Sigma_n)$ for all $n\geq 3$; for example $x_1x_2$ and
$x_2x_1$ does not commute in $\mathcal{SB}_{n}$, but their images
$\pi(x_1x_2)$ and $\pi(x_2x_1)$ commutes in $S\Sigma_n$.
\section{Simple Centralizer} \label{simple}

\begin{lem}$($\cite{gar}, Theorems~\emph{H}, \emph{K}$)$ \label{gar}
a) If  $x_i|_L\beta$ and $x_j|_L\beta$  for $\beta \in \mathcal{MB}_{n}$  and  $i\neq j$ then:\\
\emph{i)} $x_ix_j|_L\beta$ for $|i-j|\geq 2$ and\\
\emph{ii)} $x_ix_jx_i|_L\beta$ for $|i-j|=1.$

b) If  $x_i|_R\beta$ and $x_j|_R\beta$  for $\beta \in \mathcal{MB}_{n}$  and  $i\neq j$ then:\\
\emph{i)} $x_ix_j|_R\beta$ for $|i-j|\geq 2$ and\\
\emph{ii)} $x_ix_jx_i|_R\beta$ for $|i-j|=1.$\\
\end{lem}
\begin{lem}\label{gartyp} For $\beta,\gamma\in \mathcal{MB}_{n}
$ we have\\
a) if  $x_{n-1}|_L\beta\gamma$ and $x_{n-1}\nmid\beta$ then $x_{n-1}|_L\gamma,$\\
b) if  $x_{n-1}|_R\gamma\beta$ and $x_{n-1}\nmid\beta$ then $x_{n-1}|_R\gamma.$\\
\end{lem}
\proof a) The induction on the length $|\beta|$ of $\beta$ starts at
$|\beta|=1$: for $\beta=x_i$, Lemma~\ref{gar}a)i) implies
$x_i\gamma=x_ix_{n-1}\delta$ for some $\delta\in\mathcal{MB}_{n}$.
By cancellation property, $x_{n-1}|_L\gamma.$ For $|\beta|$, suppose
$\beta=x_i\beta_1$ implies $x_{n-1}|_L\beta_1\gamma.$ By
induction hypothesis we have $x_{n-1}|_L\gamma$.\\
b) The proof is symmetric to a).
\endproof

\emph{Proof of the Theorem \ref{gen}:} If $x_j\nmid\beta$ then, by
the presentation~(\ref{classical}), $\beta x_i=x_i\beta$ and
$\beta\in C_n(x_i)$. Conversely, apply induction on the length
$|\beta|$ of $\beta$. For $|\beta|=1$, $x_ix_j=x_jx_i$ implies
$|i-j|\geq 2$ or $i=j$ by the presentation~(\ref{classical}). For
$|\beta|, $ suppose on contrary, i.e., $\beta\in C_n(x_i)$ and
$x_j\mid\beta.$ Let $\beta=x_k\beta_1$ and consider the two cases:
a) $k=j$ and b) $k\neq j$. For a), by simplicity of $\beta$ we have
$x_j\nmid\beta_1x_i$ and $x_j\beta_1x_i=x_ix_j\beta_1$ but then, by
Lemma~\ref{gar}, $x_j\beta_1x_i=x_jx_ix_j\beta_2$ for some positive
braid $\beta_2$. Hence $\beta_1x_i=x_ix_j\beta_2$, which contradicts
that $x_j\nmid\beta_1x_i$. For b), $x_j\mid\beta$ and
$x_k\beta_1x_i=x_ix_k\beta_1=x_kx_i\beta_1$ implies
$\beta_1x_i=x_i\beta_1$ which contradicts the induction hypotheses.
\hfill $\Box$

\begin{lem}\label{lr} If $\beta\in \mathcal{SB}_{n}$ and $x_{n-1}|\beta$, then either $x_{n-1}|_R\beta$ or
$x_{n-1}|_L\beta$.
\end{lem}
\proof Since $\beta$ contains $x_{n-2}$ at most once so $x_{n-2}$
must be either on the right of $x_{n-1}$ or on the left of it in
$\beta$. If $x_{n-2}$ is on the right then, by the
presentation~(\ref{classical}), $x_{n-1}$ can be moved to the left
most in $\beta$. Similarly $x_{n-1}$ can be moved to the right most
in $\beta$, if  $x_{n-2}$ is on the left.
\endproof
\begin{lem}\label{consec}
If $\beta\in\mathcal{SB}_{n}$ and $\alpha\in\mathcal{SB}_{n+1}$ such
that $x_{n-1}|\beta$ and $x_{n}|\alpha$ then $\beta\alpha\neq \alpha
\beta.$
\end{lem}
\proof Suppose on contrary.
By Lemma~\ref{lr}, the following four cases have to be dealt with:\\
a) $\beta=x_{n-1}\beta_1,$ $\alpha=x_n\alpha_1$; b)
$\beta=x_{n-1}\beta_1,$ $\alpha=\alpha_2x_n$; c)
$\beta=\beta_2x_{n-1},$ $\alpha=x_n\alpha_1$; and d)
$\beta=\beta_2x_{n-1},$ $\alpha=\alpha_2x_n$, where $\beta_1,
\beta_2\in\mathcal{SB}_{n-1}$ and $\alpha_1,
\alpha_2\in\mathcal{SB}_{n}.$

For a), $x_{n-1}\beta_1x_n\alpha_1=x_n\alpha_1x_{n-1}\beta_1$. By
Theorem~\ref{gen},
$x_{n-1}x_n\beta_1\alpha_1=x_n\alpha_1x_{n-1}\beta_1$. By
Lemma~\ref{gar} and the cancellation property,
$\alpha_1x_{n-1}\beta_1=x_{n-1}x_n\delta$ for some
$\delta\in\mathcal{MB}_{n+1}$ which means $x_n|\alpha_1$ or
$x_n|\beta_1$, a contradiction.

For b),
\begin{equation}\label{2}
    x_{n-1}\beta_1\alpha_2x_n=\alpha_2x_nx_{n-1}\beta_1.
\end{equation}
By Lemma~\ref{gartyp}b, $\alpha_2x_nx_{n-1}=\delta x_n$ for some
$\delta\in\mathcal{MB}_{n+1}$. By Lemma~\ref{gar}b)ii) and the
cancellation property, $\alpha_2=\delta_1 x_{n-1}$ for some
$\delta_1\in\mathcal{MB}_{n+1}$.\linebreak Using equation
~(\ref{2}),
$x_{n-1}\beta_1\alpha_2x_n=\delta_1x_{n-1}x_nx_{n-1}\beta_1=\delta_1x_nx_{n-1}x_n\beta_1=\delta_1x_nx_{n-1}\beta_1x_n.$
We get $x_{n-1}\beta_1\alpha_2=\delta_1x_nx_{n-1}\beta_1$ which
means $x_n|\alpha_2$ or $x_n|\beta_1$, a contradiction.

 The proofs of the cases  c) and d) are symmetric to a) and b).

\endproof
The following holds in an arbitrary monoid $\mathcal{M}$ with the
cancellation law.
\begin{lem}\label{nb}
 If $\beta(\gamma_1\gamma_2)=(\gamma_1\gamma_2)\beta$ and
$\beta\gamma_1=\gamma_1\beta$ or $\beta\gamma_2=\gamma_2\beta$, then
$\beta\gamma_2=\gamma_2\beta$ or $\beta\gamma_1=\gamma_1\beta$
respectively.
\end{lem}
\begin{lem}\label{nm}
If $\beta\in\mathcal{SB}_{n}$ and $x_{n-1}|\beta,$ then
$x_{n}\nmid\gamma$ for any $\gamma\in C_{n+m}(\beta)$, where $m\geq
1.$
\end{lem}
\proof By Lemma~\ref{lr} all the divisors $x_j$, $j\geq n+1$ can be
written at the beginning or at the end, i.e, $\gamma=\eta\rho\mu$
such that $\rho\in\mathcal{SB}_{n+1}$, $x_j\nmid\eta$ and
$x_j\nmid\mu$ for $j\leq n.$ Since $\beta\eta=\eta\beta$ and
$\beta\mu=\mu\beta$, so by Lemma~\ref{nb}, we find
$\beta\rho=\rho\beta$. By Lemma~\ref{consec}, $x_n\nmid\rho.$
\endproof
 \emph{Proof of the Theorem \ref{struc}:}
Follows from Lemma~\ref{nm} and the presentation~(\ref{classical}).

\begin{lem}\label{x_{n-1}}
a) $c_2(x_1)=2$ and\\
 b) $c_{n}(x_{n-1})=2F_{2n-5}$ $\forall\,\, n\geq 3 $.
\end{lem}
\proof a) Clearly $C_n(x_1)=\mathcal{SB}_{2}=\{e,x_1\}$.

b) Since $\mathcal{SB}_{n-2}. x_{n-1}= x_{n-1}.\mathcal{SB}_{n-2}$,
so by Theorem~\ref{gen},
$$C_{n}(x_{n-1})=\mathcal{SB}_{n-2}\amalg
x_{n-1}.\mathcal{SB}_{n-2}.$$
\endproof
 \emph{Proof of the Proposition\ref{c(x)}:}
 By Corollary~\ref{cent} and Lemma~\ref{x_{n-1}}. \hfill $\Box$

\section{Commuting  graphs} \label{commuting}
Clearly the degree of $\beta\in \mathcal{SB}_{n}$ in
$\Gamma(\mathcal{SB}_{n})$ is equal to $c_n(\beta)-2$: we minus $e$
and $\beta$ from $C_n(\beta)$ in order to make the vertices
non-central and avoid loop in $\Gamma(\mathcal{SB}_{n}).$ The graph
$\Gamma(\mathcal{SB}_{n-1})$ has a canonical embedding in
$\Gamma(\mathcal{SB}_{n}).$ Similarly, The graphs
$\Gamma(S\Sigma_{n-1})$ and $\Gamma(\Sigma_{n-1})$ are embedded in
$\Gamma(S\Sigma_n)$ and $\Gamma(\Sigma_n)$ respectively.

The canonical form for a positive braid (the smallest element in the
equivalence class of a braid) was introduced in~\cite{bok, barbu} to
solve the word problem in $\mathcal{MB}_{n}.$ In order to generate
the $\Gamma(\mathcal{SB}_{5})$, we used these canonical forms,
Theorem~\ref{gen}, Lemma~\ref{consec}, and Lemma~\ref{nb}.
\begin{examples}
1) $\alpha=x_1x_3x_2$ and $\beta=x_2x_3$ does not commute because
the canonical form of $\alpha\beta$ and $\beta\alpha$ is
$x_1x_3x^2_2x_3$ and $x_2x_1x^2_3x_2$ respectively.

2) $\alpha=x_1x_3$ and $\beta=x_2x_4x_3$ does not commute by
Lemma~\ref{consec}.

3) $\alpha=x_1x_4$ and $\beta=x_2x_1$ does not commute by
Lemma~\ref{nb}.
\end{examples}
\begin{lem}$($\cite{doug}, page 250, example \emph{7.1.5.}$)$
\label{k33} The complete bipartite graph $K_{3,3}$ and the complete
graph $K_5$ are non-planar.
\end{lem}
\emph{Proof of the Proposition~\ref{planner}:} a) The graph
$\Gamma(\mathcal{SB}_{6})$ is non-planar (by Lemma~\ref{k33})
because it contains $K_{3,3}$ as a subgraph:
\begin{center}
\protect
\begin{picture}(290,30)%%    fig 5b  %%
\put(106,20){$x_1$}\put(100,-10){$x_1x_2$} \put(100,-40){$x_2x_1$}
\put(207,20){$x_4$}\put(201,-10){$x_4x_5$} \put(201,-40){$x_5x_4$}
%\put(45,30){$K_{3,3}:$}
\put(127.3,18){$\bullet$}\put(125,-13){$\bullet$}\put(125,-45){$\bullet$}
\put(192,18){$\bullet$}\put(192,-14){$\bullet$}\put(192,-45){$\bullet$}
\put(128,-44){\line(1,1){65}} \put(193,-42){\line(-1,1){62}}
\put(129,-11){\line(1,0){64}}
\put(128,21){\line(1,0){64}}\put(128,-42){\line(1,0){64}}
\put(128,-43){\line(2,1){65}} \put(128,21.5){\line(2,-1){65}}
\put(128,-10){\line(2,-1){65}}
 \put(129,-11){\line(2,1){65}}
\end{picture}\protect
\end{center}
\vspace{2.5cm}
 There are $33$ non-central elements in
$\mathcal{SB}_{5}.$ The connected component of
$\Gamma(\mathcal{SB}_{5})$ is given in the following diagram and the
rest of the vertices are isolated
\begin{figure}[h]
\begin{center}
\includegraphics[width=6.5cm]{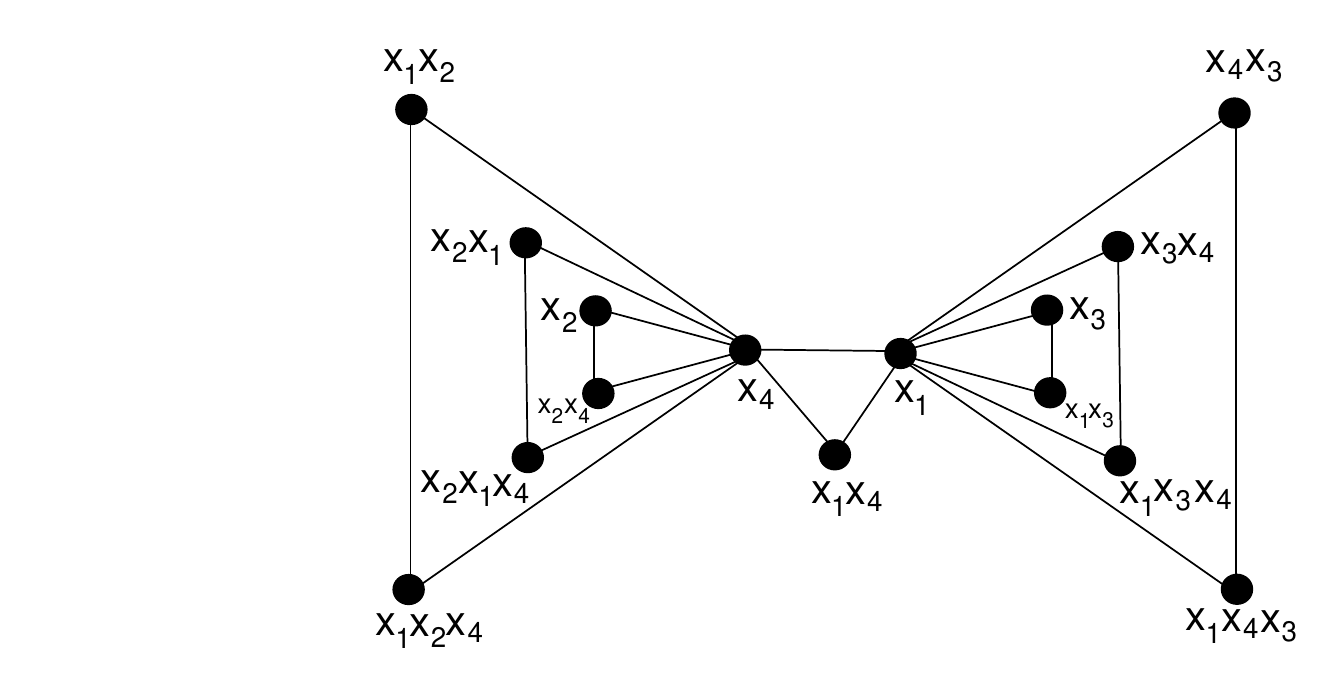}\\
%\caption{$C_{\Gamma}$}\label{1}
\end{center}
\end{figure}

b) The graph $\Gamma(S\Sigma_5)$ is non-planar (by Lemma~\ref{k33})
because it contains $K_5$ as a subgraph:
\begin{figure}[h]
\begin{center}
\includegraphics[width=4cm]{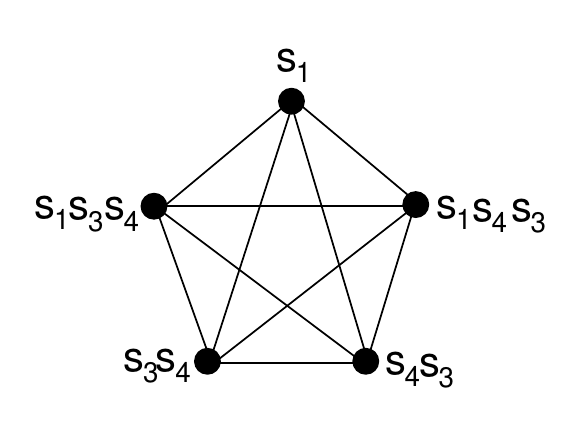}\\
%\caption{$K_{5}$}\label{1}
\end{center}
\end{figure}

 It is easy to check that $\Gamma(\Sigma_4)$ is planar:
$\Gamma(\Sigma_4)$ contains the following cycle and the rest of
vertices have degrees at most one.
\begin{figure}[h]
\begin{center}
\includegraphics[width=4cm]{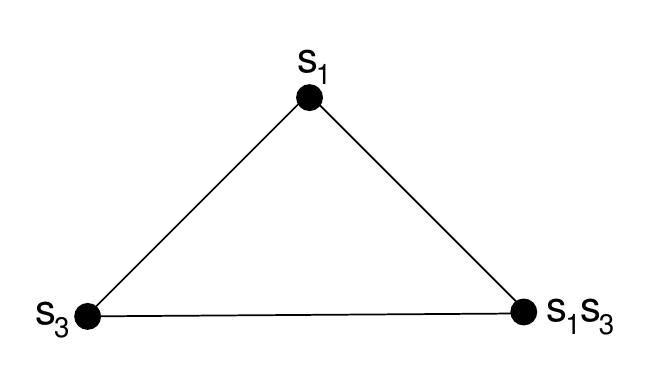}\\
%\caption{$L$}\label{1}
\end{center}
\end{figure}
\endproof
%\vspace{2cm}
 \hfill $\Box$\\
 %\textbf{Acknowledgement:} The authors are thankful to the referee for his/her valuable
% suggestions to improve
% the presentation of the paper.

\end{document}